\documentclass[graybox]{svmult}

\usepackage{type1cm}        % activate if the above 3 fonts are
\usepackage{makeidx}         % allows index generation
\usepackage{graphicx}        % standard LaTeX graphics tool
\usepackage{multicol}        % used for the two-column index
\usepackage[bottom]{footmisc}% places footnotes at page bottom

\usepackage{newtxtext}       % 
\usepackage[varvw]{newtxmath}       % selects Times Roman as basic font

\usepackage{lineno,hyperref}
\modulolinenumbers[5]

\usepackage{multirow}
\usepackage{slashbox}
\usepackage{subfigure}

\newcommand\e[1]{{\rm e}^{#1}}
\renewcommand{\E}{\mathbb{E}}

\makeindex             % used for the subject index
\usepackage{dsfont}

\usepackage{xcolor}

\begin{document}

\title*{IMEX-RK finite volume methods for nonlinear 1d parabolic PDEs. Application to option pricing}
\titlerunning{IMEX-RK FV methods for nonlinear 1d parabolic PDEs in finance}
\author{J. G. L\'opez-Salas 
        and M. Su\'arez-Taboada
        and M. J. Castro
        and A. M. Ferreiro-Ferreiro
        and J. A. Garc\'ia-Rodr\'iguez
        }
\authorrunning{JG L\'opez-Salas, M Su\'arez-Taboada, A Ferreiro-Ferreiro, JA Garc\'ia, MJ Castro}

% \author[auth1]{J. G. L\'opez-Salas}\ead{jose.lsalas@udc.es}
% \address[auth1]{Department of Mathematics, Faculty of Informatics and CITIC, Campus Elvi\~na s/n, 15071-A Coru\~na (Spain)}

% \author[auth1]{M. Su\'arez-Taboada}\ead{maria.suarez3@udc.es}

% \author[auth2]{M. J. Castro\corref{cor1}}\ead{mjcastro@uma.es}
% \cortext[cor1]{Corresponding author.}
% \address[auth2]{Department of An\'alisis Matem\'atico, Facultad de Ciencias, University of M\'alaga, Campus de Teatinos s/n, M\'alaga, 29080-Andaluc\'ia (Spain)}

% \author[auth1]{A. M. Ferreiro-Ferreiro}\ead{ana.fferreiro@udc.es}
% \author[auth1]{J. A. Garc\'ia-Rodr\'iguez}\ead{jose.garcia.rodriguez@udc.es}

\institute{J. G. L\'opez-Salas, M. Su\'arez-Taboada, A. M. Ferreiro-Ferreiro, J. A. Garc\'ia-Rodr\'iguez \at Department of Mathematics, Faculty of Informatics and CITIC, Campus Elvi\~na s/n, 15071-A Coru\~na (Spain), \email{jose.lsalas@udc.es, maria.suarez3@udc.es, ana.fferreiro@udc.es, jose.garcia.rodriguez@udc.es}
\and M. J. Castro \at Department of An\'alisis Matem\'atico, Facultad de Ciencias, University of M\'alaga, Campus de Teatinos s/n, M\'alaga, 29080-Andaluc\'ia (Spain)  \email{mjcastro@uma.es}}

\maketitle
\vspace*{-2.5cm}
\abstract*{   
    The goal of this paper is to develop 2nd order Implicit-Explicit Runge-Kutta (IMEX-RK)  finite volume (FV) schemes for solving 1d parabolic PDEs for option pricing, with possible nonlinearities in the source and advection terms. The spatial semi-discretization of the advection is carried out by combining finite volume methods with 2nd order  state reconstructions; while the diffusive terms are discretized using second order finite differences.  The time integration is performed by means of IMEX-RK time integrators: the advection is treated explicitly, and the diffusion, implicitly.  The obtained numerical schemes have several advantages: they are computationally very efficient, thanks to the implicit discretization of the diffusion in the IMEX-RK time integrators, that allows to overcome the strict time step restriction;  they yield second order accuracy for even nonlinear problems and with non-regular initial conditions; and they can be extended to higher order.
}
\abstract{
    The goal of this paper is to develop 2nd order Implicit-Explicit Runge-Kutta (IMEX-RK)  finite volume (FV) schemes for solving 1d parabolic PDEs for option pricing, with possible nonlinearities in the source and advection terms. The spatial semi-discretization of the advection is carried out by combining finite volume methods with 2nd order  state reconstructions; while the diffusive terms are discretized using second order finite differences.  The time integration is performed by means of IMEX-RK time integrators: the advection is treated explicitly, and the diffusion, implicitly.  The obtained numerical schemes have several advantages: they are computationally very efficient, thanks to the implicit discretization of the diffusion in the IMEX-RK time integrators, that allows to overcome the strict time step restriction;  they yield second order accuracy for even nonlinear problems and with non-regular initial conditions; and they can be extended to higher order.
    % Finally, it is straightforward to extend the presented numerical schemes to higher order, as they are based in well studied and solid numerical techniques, like IMEX-RK time integrators, numerical fluxes and state reconstruction operators.
}

\vspace{-0.5cm}
\section{Introduction}

The goal of this paper is to present a general technique for building  high-order FV numerical schemes for solving 1d parabolic PDEs arising in option pricing. The advantages of the here proposed schemes are:
they are high order, they are computationally very efficient and they are very general and extendable.
The first important novelty is that the proposed numerical schemes are able to preserve high-order, even when the initial condition is not smooth. And even more important, as fixed point schemes are no longer needed to deal with the nonlinearities, the schemes preserve high-order when we have nonlinearities in the advection and/or the source term.
Additionally, the obtained schemes only need to satisfy the convection stability condition. Therefore, larger time steps can be taken, thus the here proposed schemes being orders of magnitude more efficient than explicit schemes, from the computational point of view.
Finally, the here proposed numerical schemes rely on very well established finite volume tools for dealing with the convective terms at high-order: upwind of the convection is treated by means of numerical fluxes based on Riemann solvers, and well known state reconstruction operators are used for increasing the order of convergence.
To the best of our knowledge this is the first time this systematic framework is presented for building higher order schemes for solving nonlinear parabolic PDEs in option pricing.

More precisely, the focus of this article is to develop numerical schemes that are able to solve very efficiently nonlinear PDEs in finance, specifically semilinear PDEs and second-order quasilinear PDEs. These kinds of PDEs are of great relevance  in mathematical finance, and even more in the post financial 2008 crisis scenario, where the computation of value at risk has become a key task. 
These models are given by advection-diffusion-reaction 1d scalar PDEs, under the general expression:
\begin{equation}
\label{AdvDiffPDE}
\dfrac{\partial u}{\partial t}+
a(x,t)\dfrac{\partial u}{\partial x}
%----------------
+b(x,t)\dfrac{\partial^2 u}{\partial x^2}
%------------------
+H\left(x,t,u,\dfrac{\partial u}{\partial x}\right)=0,
\end{equation}
with $H$ a continuous function.

The discretization of these kind of financial PDEs, for the linear case, using finite difference and finite element methods, has been studied in  \cite{Duffy06,Wilmott06,Pironneau05,DuffyExpFitted13} and references therein. 
% The combination of finite differences with Exponentially Fitted techniques is explained in \cite{DuffyExpFitted13}. Besides, Alternate Directions (ADI) with finite differences is illustrated in \cite{adi}.
% \cite{Foulon10}. 
However, the development of finite difference and finite element numerical methods for the PDEs arising in mathematical finance presents several well known difficulties. First, and most important, these numerical methods usually show instabilities when the advection term becomes larger and/or the diffusion operator is degenerate. Upwinding techniques are needed to overcome this issue. Secondly, developing high order pricers is challenging, because 2nd order (or higher) convergence is lost when the initial condition is not regular: this is the typical situation in option pricing, as the initial condition is given by a payoff that is usually singular. 
% Finally, another difficulty, derived from the previous ones, is achieving accurate and non oscillatory approximations of the Greeks. The derivatives of the solution are usually computed by means of finite difference formulas, which are very sensitive to small errors in the approximation of the prices. The higher the derivatives the more difficult is obtaining approximations without oscillations. This question is of paramount importance, since the Greeks are vital for trading purposes. Developing very accurate and high order schemes is a key step towards attaining non-fluctuating approximations of the Greeks.

% In order to avoid the problems originated by non-smooth payoffs, smoothing techniques working on irregular initial data were proposed in the literature, see \cite{WADE2007144}. A commonly used smoothing technique is the so-called Rannacher's method, see \cite{Rannacher}. It is well known that the second order Crank-Nicolson time marching scheme loses order when initial conditions are non-smooth, or the initial and boundary conditions are discontinuous, which is the situation with barrier options. Rannacher proposed a way to suppress wrongful initial oscillations, by preceding Crank-Nicolson with a few implicit steps.
Several numerical approaches have been presented in the literature in order to deal with convection  in the context of finance PDEs. 
% One approach is the method of characteristics. In \cite{ForsythCharact06}, Forsyth et al. solve option pricing problems with finite differences combined with the semi-Lagrangian characteristics method. 
% In the finite element setup, semi-Lagrangian characteristics was applied in \cite{ARREGUI2018}
% \cite{VazquezNogueiras2006-1,VazquezNogueiras2006-2} 
% for pricing XVA.
% In \cite{Ferreiro13} the authors present a semi-Lagrangian finite difference method for pricing business companies. 
% The main disadvantages of semi-Lagrangian methods in option pricing is the difficulty to build high order numerical schemes. In fact, these numerical methods do not achieve second order convergence in the presence of non-smoothness of either the payoff or the boundary conditions, or in the nonlinear case. On top of that, the computational cost of characteristic method is high due to the demanding compulsory search at the foot of the characteristic and the required interpolation. 
One technique for the upwinding of the advection terms, is the use of finite volume methods. One of the first works in  applying finite volume methods for option pricing was \cite{ForsythZvan01}. 
% Later, in \cite{mathias13}  conservative explicit finite volume methods were proposed for convection dominated pricing problems. More precisely, the authors propose to use the extension of the central schemes presented by Nessyahu-Tadmor in \cite{NT1990} to the advection-diffusion problem developed in \cite{KT2000}.
%In \cite{mathias13} the KT explicit finite volume method is applied to the PDE equations in conservative form.
Recently, in \cite{Tadmor2ndOrder2019}, conservative explicit finite volume schemes were presented for  dealing with convection in the context of pricing problems. More precisely, the authors propose a 
second order method with appropriate time methods and slope limiters, that is based on the central schemes presented by Nessyahu-Tadmor in \cite{NT1990} and \cite{KT2000}.
% The authors propose a 
% second order method with appropriate time methods and slope limiters. 
In \cite{Tadmor2019} the authors apply the explicit third order Kurganov-Levy scheme presented in \cite{KL2000} along with the CWENO reconstructions presented in \cite{LPR1999}. In all theses articles, it is shown that explicit finite volume schemes do not suffer loss in the order of convergence. Besides, they are able to obtain approximations of the Greeks, without oscillations. In finance, the Greeks are the sensitivities of the option price with respect to changes in its underlying parameters.
Nevertheless, these works present numerical schemes that are explicit in time. Explicit time integrators introduce a severe restriction in the time step, imposed by the Von Neuman stability condition related to the diffusion terms. As a consequence, these schemes have a huge computational cost for refined meshes in space.
% and are not able in practice to work with refined meshes in space, specially in problems with spatial dimension greater than one.

In this work we develop general finite volume numerical solvers for option pricing problems. The proposed schemes address the upwinding problems of finite difference and finite element methods, while at the same time allowing to use large time  steps in the time discretization. More precisely, we present a general technique, following \cite{Russo05,boscarino,russo-finance}, for building  second-order Implicit-Explicit (IMEX) Runge-Kutta finite volume solvers for option pricing. This numerical scheme allows to use different numerical flux functions and opens the door to the development of higher order schemes based in state reconstructions, in mathematical finance.  Also the proposed method is able to overcome the severe time step restriction thanks to the implicit treatment of the diffusive part, while retaining at the same time  the benefits of treating the advective term by means of an explicit finite volume scheme. In this way the stability condition of the IMEX scheme allows to use the same time step of the advective part, which is far larger than the diffusive time step. Moreover, finite volume schemes allow to address the loss of order of convergence when initial data is non-smooth, since they handle the integral version of the equations, working with the averaged solutions in each cell. Consequently, true second order schemes are proposed for option pricing problems, that also allow to recover accurate and non oscillatory approximations of the Greeks. Besides, we show that IMEX schemes are particularly well suited for solving PDEs with nonlinearities in the convection and/or source terms, since they treat these terms explicitly while use implicit discretization only for the linear diffusion terms. This approach avoids either inefficient fixed point iterations or complex and less accurate linearization algorithms. 

\vspace{-0.5cm}
\section{Numerical methods. Finite volume  IMEX Runge-Kutta\label{sec:numerical-schemes}}
In this section we present a second order finite volume semi-implicit numerical scheme for solving \eqref{AdvDiffPDE}.
First, if possible, the equation \eqref{AdvDiffPDE} must be written in conservative form:
\begin{equation} \label{ConsForm1D}
    \dfrac{\partial u}{\partial t}  + \dfrac{\partial}{\partial s} f(u,s)=\dfrac{\partial}{\partial s} g(u_s,s)+h(u).
\end{equation}
If the PDE \eqref{AdvDiffPDE} cannot be written in conservative form, special numerical methods based in \cite{castro2017wella} for treating the nonlinearity can also be developed, taking also advantage of all the techniques discussed in the present paper.
% taking $x\equiv s$.
The numerical solution of equation \eqref{ConsForm1D} using a explicit finite volume scheme may have a huge computational cost because of the tiny time steps induced by the diffusive terms. To avoid this difficulty we consider IMEX Runge-Kutta methods (see \cite{Russo05}). 
% These methods play a major rule in the treatment of differential systems governed by stiff and non stiff terms. 
% of the form:
% \begin{equation}\label{eq:stiffedo}
% \dfrac{\partial U}{\partial t}  + F (U)= S  (U),
% \end{equation}
% where $F,S:R^N\to R^N$, being $F$  the non-stiff term and $S$  the stiff one.
% % This equation involves an explicit and a implicit part.
% Both can be solved simultaneously using the same {\em time step} of the convective part,
% by means of a second order IMEX Runge-Kutta  scheme.
% In this work we apply this technique to  equation \eqref{ConsForm1D} following \cite{Russo05}. 
The procedure for obtaining the numerical scheme can be summarized as follows. First, we perform a spatial finite volume semi-discretization of \eqref{ConsForm1D}, explicit in convection and reaction, and implicit in the diffusive part. As a result we obtain a stiff time ODE system, that we discretize using IMEX Runge-Kutta methods.  
% In what follows we succinctly describe the space and time discretizations.

\subsection{Finite volume spatial semi-discretization}

The spatial semi-discreti\-zation of the advective and source terms is performed by means of an explicit finite vo\-lume scheme. The spatial domain is split into cells $\{I_i\}$,  with $I_i=[s_{i-1/2},s_{i+1/2}]$, $i=\ldots,-1,0,1, \ldots$, being $s_i$ the center of the cell $I_i$. Let $|I_i|$ be the size of cell $I_i$.
The basic unknowns of our problem are the averages of the solution $u(s,t)$ in the cells $\{I_i\}$, $\bar{u}_i= \frac{1}{|I_i|}\int_{I_i} u \, ds$.
Integrating equation \eqref{ConsForm1D} in space on $I_i$ and dividing by $|I_i|$ we obtain the semi-discrete finite volume scheme:
{\small
\begin{align}
\dfrac{d \bar{u}_i}{dt}
 =&-\dfrac{1}{|I_i|} \left( \mathcal{F}_{i+1/2}(t) - \mathcal{F}_{i-1/2}(t)\right) 
 \label{eq:convectiveterms}
\\
&+\dfrac{1}{ |I_i|} \left( \mathcal{G}_{i+1/2}(t)  - \mathcal{G}_{i-1/2}(t)\right)
\label{eq:difussionterms}
\\
&+\dfrac{1}{ |I_i|} \int_{I_i} h(\mathcal{R}^t(s)) \, ds.  \label{eq:sourceterms}
\end{align}
}
In the previous expression, $\mathcal{F}_{i+1/2}$(t) is a numerical flux evaluated on the reconstructed states $u^\pm_{i+1/2}(t)$ at the intercell $s_{i+1/2}$. That is
$u^\pm_{i \pm1/2}(t)= \lim_{s\to s^\pm_{i \pm 1/2}} \mathcal{R}^t(s), $
where $\mathcal{R}^t(s)$ is a reconstruction operator of the  unknown function $u(s, t)$ defined in terms of the cell averages $\{\bar{u}_i(t)\}$, that could be written generally as
$\mathcal{R}^t(s) = \sum_i P^t_i (s) \mathds{1}_{s\in I_i},$
where $P^t_i(s)$ is a reconstruction polynomial at cell $I_i$ at time $t$ satisfying some accuracy and non oscillatory properties, defined from the cell averages $\{\bar u_j(t)\}$ on a given stencil $j \in S_i$, and $\mathds{1}_{s\in I_i}$ is the indicator function of cell $I_i$. For second order schemes, the reconstruction operator has to be at least piecewise linear. 
For example the left reconstructed value at the edge $s_{i+1/2}$ is
$
u^-_{i + 1/2}= \bar{u}_{i}+\frac{1}{2} u'_{i},
$
where the slope $u'_{i}$ is a first order approximation of the space derivative of $u(s,t)$ at point $s_i$ at every time $t$.
Slope limiters must be applied to satisfy the TVD property: we use the minmod limiter.
% , where the slope is given by
% $$
% u'_{i}= {\rm minmod} ( \bar{u}_{i}-\bar{u}_{i-1},\bar{u}_{i+1}-\bar{u}_{i}),
% $$
% $$
% {\rm minmod}(a,b)=
% \begin{cases}
% \min (a,b) & \text{if } a,b>0,\\
% \max(a,b) & \text{if } a,b<0,\\
% 0 & \text{otherwise\,. } 
% \end{cases}
% $$
In this work we consider a LLF numerical flux:
{\small
$$
\mathcal{F}_{i+1/2}=\dfrac{1}{2}\left(f(u^-_{i+1/2},s_{i+1/2})+f(u^+_{i+1/2},s_{i+1/2})\right)- \dfrac{\alpha_{i+1/2}}{2}\left(u^+_{i+1/2}-u^-_{i+1/2}\right),
$$
}
with
$
\alpha_{i+1/2}=\left|\frac{\partial f}{\partial u}\left(\frac{u^-_{i+1/2}+u^+_{i+1/2}}{2}\right)\right|.
$
% In the previous expression we have dropped the time dependency for simplicity. 
The integral of the source term \eqref{eq:sourceterms} can be explicitly discretized using a quadrature rule. In this work, as we are interested on second order schemes we use the midpoint rule
$
\int_{I_i}  h(\mathcal{R}^t(s)) ds\approx |I_i|\, h(\mathcal{R}^t(s_i))=|I_i|\,h(\bar{u}_i).
$
Finally, $\mathcal{G}_{i+1/2}(t)$ is defined by
$
\mathcal{G}_{i+1/2}=g\left( \frac{\bar{u}_{i+1}- \bar{u}_{i}}{|I_i|},s_{i+1/2}\right).
$
Again, we have dropped the time dependency for simplicity. Observe that $\frac{1}{|I_i|}\left(\mathcal{G}_{i+1/2}-\mathcal{G}_{i-1/2}\right)$
is a second order approximation of the diffusive term $\frac{\partial}{\partial s} g(u_s,s)$ at $s=s_i$.

\vspace{-0.75cm}
\subsection{IMEX Runge-Kutta time discretization}
After performing the spatial semi-discretization of equation \eqref{ConsForm1D} we obtain a stiff ODE system of the form
\begin{equation}\label{eq:stiffedo}
\dfrac{\partial U}{\partial t}  + F (U)= S  (U),
\end{equation}
where \mbox{$U=(\bar{u}_i(t))$} and $F,S:\mathbb{R}^N\to \mathbb{R}^N$, being $F$ the non-stiff term and $S$  the stiff one.
An IMEX scheme consists of applying an implicit discretization to the stiff term
and an explicit one to the non stiff term. In this way, both can be solved simultaneously with high order accuracy using the same {\em time step} of the convective part, which is in general much larger than the time step of the diffusive part.
When IMEX is applied to system \eqref{eq:stiffedo} it takes the form:
{\small
\begin{align}
U^{(k)} &= U^n -\Delta t \sum_{l=1}^{k-1} \tilde{a}_{kl} F(t_n+\tilde{c}_l\Delta t,U^{(l)})+\Delta t \sum_{l=1}^{\rho} a_{kl} S(t_n+{c}_l\Delta t,U^{(l)}), \\
U^{n+1} &= U^n- \Delta t \sum_{k=1}^{\rho} \tilde{\omega}_{k} F(t_n+\tilde{c}_k\Delta t,U^{(k)})+\Delta t \sum_{k=1}^{\rho} \omega_{k}S(t_n+{c}_k\Delta t,U^{(k)}),
\end{align}
}
where $U^n=(\bar{u}_i^n)$, $U^{n+1}=(\bar{u}_i^{n+1})$ are the vectors of the unknown cell averages at times $t^n$ and $t^{n+1}$, thus $U^{(k)}$ and $U^{(l)}$ are the vector of unknowns at the stages $k,l$ of the IMEX method. The matrices $\tilde{A} = (\tilde{a}_{kl})$, with $\tilde{a}_{kl} = 0$ for $l \geq k$, and $A=(a_{kl})$ are square matrices of order $\rho$
% such that the ensuing scheme is implicit in $S$ and explicit in $F$. Solving efficiently at each time step the system of equations corresponding to the implicit part is extremely important. Therefore, one usually considers $a_{kl} = 0$, for $l > k$, the so-called diagonally implicit Runge-Kutta (DIRK) schemes . 
that are given by a double Butcher tableau:
\begin{center}
\begin{tabular}{c|c}
$\tilde{c}$  & $\tilde{A}$ \\
\hline
&  $\tilde{\omega}$ 
\end{tabular},
\quad
\begin{tabular}{c|c}
$c$  & $A$ \\
\hline
         &  $\omega$ 
\end{tabular},
\end{center}
% where  $\tilde{w}=(\tilde{w}_1,\ldots,\tilde{w}_{\rho})$ and $w=(w_1,\ldots,w_{\rho})$. Besides, the coefficient vectors $\tilde{c}=(\tilde{c}_1,\ldots,\tilde{c}_{\rho})^T$ and $c=({c}_1,\ldots,{c}_{\rho})^T$ are only used for the treatment of non autonomous systems, and have to satisfy the relations
% \begin{equation}
%   \tilde{c}_k=  \sum_{l=1}^{k-1} \tilde{a}_{kl}, \quad c_k=\sum_{l=1}^{k} {a}_{kl}.
% \end{equation}
In this work we will consider the second order IMEX-SSP2(2,2,2) L-stable scheme (see \cite{Russo05}) that is given by
\begin{center}
\begin{tabular}{c|cc}
0  & 0 & 0 \\
1  & 1 & 0 \\
\hline
& $1/2$ & $1/2$  
\end{tabular}
\quad
\begin{tabular}{c|cc}
$\gamma$  & $\gamma$ & 0 \\
$1-\gamma$  & $1-2\gamma$ & $\gamma$ \\
\hline
& $1/2$ & $1/2$  
\end{tabular}
\quad 
$\gamma=1-\dfrac{1}{\sqrt{2}}.$ 
\end{center}
An explicit  time integrator needs extremely small time steps due to the following stability conditions
$
    \eta \frac{\Delta t}{ (\Delta s)^2} \leq \frac{1}{2}, \label{eq:stabilityDiffusion}
$
and
$
    \alpha \frac{\Delta t}{\Delta s} \leq 1\, , \label{eq:stabilityConvection}
$ 
where  $\eta=\left\lvert\frac{\partial g}{\partial u_s}\right\rvert,\, \alpha = \left\lvert\frac{\partial f}{\partial u} \right\rvert ,$ for all cells $I_{i}$ and for all boundary points $s_{i\pm 1/2}$.
However, IMEX only needs to satisfy the advection stability condition \eqref{eq:stabilityConvection}.

\vspace{-0.25cm}
\begin{figure}[!h]
\centering \includegraphics[width=0.8\linewidth]{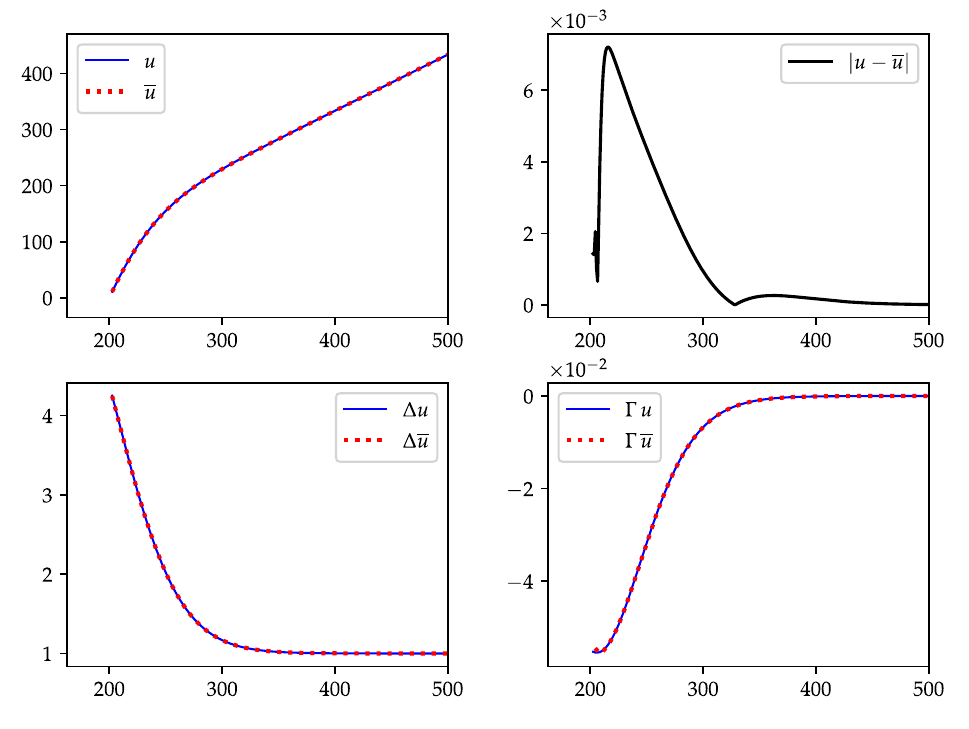}  
\caption{Down-and-out call option prices, numerical errors and Greeks ($\Delta$, $\Gamma$) at $t=T$.}
\label{fig:testBarrier}
\end{figure}

\vspace{-0.25cm}
\begin{table}[!h]
\begin{center}
{\scriptsize
\begin{tabular}{|c||c|c|c|c|}
\hline
& \multicolumn{4}{c|}{IMEX}\\
\hline
$N$  & $L_1$ \text{ error} &  \text{ Order}& $\Delta t$ & Time (s)  \\
\hline
$50$ &  $1.3889\times 10^{2}$ & $--$ &  $1.00\times 10^{0}$ & $1.8\times 10^{-4}$ \\
$100$ &  $3.4052\times 10^{1}$ & $2.03$&  $5.00\times 10^{-1}$ & $2.6\times 10^{-4}$\\
$200$ &  $8.5310\times 10^{0}$ & $2.03$ & $2.50\times 10^{-1}$ & $4.6\times 10^{-4}$ \\
$400$ &  $2.1249\times 10^{0}$ & $2.02$ & $1.25\times 10^{-1}$ &  $9.4\times 10^{-4}$ \\
$800$ &  $5.2912\times 10^{-1}$ & $2.01$ & $6.25\times 10^{-2}$ & $2.4\times 10^{-3}$ \\
$1600$ &  $1.3097\times 10^{-1}$ & $1.98$ & $3.13\times 10^{-2}$ & $7.3\times 10^{-3}$ \\
$3200$ &  $3.1547\times 10^{-2}$ & $2.00$ & $1.56\times 10^{-2}$ & $1.6\times 10^{-2}$ \\
$6400$ &  $6.7624\times 10^{-3}$ & $2.26$ &  $7.81\times 10^{-3}$ & $4.5\times 10^{-2}$ \\
\hline
\end{tabular} 
\begin{tabular}{|c||c|c|c|c|}
    \hline
& \multicolumn{4}{c|}{Explicit}\\
\hline
$N$  & $L_1$ \text{ error} &  \text{ Order}& $\Delta t$ & Time (s)  \\
\hline
$50$ & $1.3979\times 10^{2}$ & $--$  &  $5.00\times 10^{-3}$ & $1.8\times 10^{-3}$ \\
$100$ & $3.4401\times 10^{1}$ &$2.02$ &  $1.25\times 10^{-3}$ & $7.7\times 10^{-3}$\\
$200$ & $8.5373\times 10^{0}$ &$2.01$ & $3.12\times 10^{-4}$ & $3.0\times 10^{-2}$\\
$400$ & $2.1271\times  10^{0}$ &$2.01$ & $7.81\times 10^{-5}$ & $1.3\times 10^{-1}$\\
$800$ & $5.3130\times 10^{-1}$ &$2.01$ & $1.95\times 10^{-5}$ & $9.5\times 10^{-1}$\\
$1600$ & $1.3316\times 10^{-1}$ &$2.01$ & $4.88\times 10^{-6}$ & $6.4\times 10^{0}$\\
$3200$ & $3.3721\times 10^{-2}$ &$2.05$ & $1.22\times 10^{-6}$& $6.4\times 10^{1}$\\
$6400$ & $8.8809\times 10^{-3}$ &$2.22$  & $3.05\times 10^{-7}$ & $5.5\times 10^{2}$\\
\hline
\end{tabular} 
\caption{$L_1$ errors and $L_1$ orders of convergence of the IMEX and explicit finite volume methods for the down-and-out call option.}
\label{tb-Barrier} 
} 
\end{center}
\end{table}

\vspace{-1cm}
\section{Numerical experiments\label{sec:numerical-results}}
\vspace{-0.25cm}
In this section the developed numerical method is applied to the discretization and solution of linear and nonlinear financial PDEs, namely barrier options and value at risk.
% More precisely, experiments under the linear and nonlinear Black-Scholes models for vanilla, butterfly and barrier options are presented in Section \ref{sec:numerical-results-BS}. 
To show  the efficiency and accuracy, the numerical results are compared with the  analytical solutions. 
% Later, in Section \ref{sec:NumExpAsian} two dimensional problems in space are solved. Indeed, Asian options are priced. 
At each one of the following subsections, we start by describing the model and writing the involved PDE in conservative form. Then, graphs containing numerical results, such as option prices, Greeks (Delta and Gamma) and numerical errors are presented. Moreover, tables for the $L_1$ errors and the $L_1$ orders of convergence are shown. Additionally, a comparison of the time step sizes supplied by the stability conditions of the explicit and IMEX Runge-Kutta methods is presented. For all the tests a CFL of $0.5$ is considered for the stability condition.

\vspace{-0.75cm}
\subsection{Linear example. Options under the Black-Scholes model\label{sec:numerical-results-BS}}
\vspace{-0.0cm}
As a lineal example we consider a exotic barrier option.
Due to the sharp discontinuity arising at the barrier this option is mathematically interesting in the PDE world. We will price this product with our proposed finite volume IMEX Runge-Kutta schemes.

Barrier options are exotic path-dependent options. One example of barrier options is the down-and-out call option. This derivative pays $\max(s-K,0)$ at expiry, unless at any previous time the underlying asset touched or crossed a prespecified level $B$, called the barrier. In that situation the option becomes worthless. There are also \textit{in} options which only pay off if the asset reached or crossed the barrier, otherwise they expire worthless. 
% These barrier options are called continuously monitored barrier options. 
Barrier options can be modelled by the classical Black-Scholes model presented in  \cite{Merton73}-\cite{BS73},
% This is a classical model, but continues to be the gold pattern in the financial industry, and sereral kinds of options like vanilla, butterfly and exotic barrier options can be priced using this model, and fixing different payoff functions.
% In this work we  focus on pricing  European-style options: vanilla options and exotic barrier options, as well as computing the counterparty risk. 
where the dynamics of the price of the underlying asset is modeled by means of a Stochastic Differential Equation (SDE),
$\frac{ds_t}{s_t}=(r-q) dt+\sigma dW_t,
$
being $W_t$ a standard Brownian motion, the parameter $r\in \mathbb{R}$ is the risk free constant interest rate, $q\in \mathbb{R}$ is the continuous dividend yield,
% , since the coefficient on $dt$ in \eqref{eq:SDE_BS}, the so-called  mean rate of return, is considered as $r-q$. 
and the parameter $\sigma\in \mathbb{R}^+$ is the volatility of the stock price, which is again considered as constant.  This SDE describes the risk-neutral dynamics of the underlying asset price.
Applying Itô's lemma a linear parabolic backward  PDE is obtained, and after 
% Hereafter, in this work we will work forward in time by 
making a time reversal change of variable the forward in time Black-Scholes PDE is:
{\small
\begin{equation} \label{eq:ForwardBS_PDE}
    \dfrac{\partial u}{\partial t}  - \dfrac{1}{2}\sigma^2 s^2\dfrac{\partial^2 u}{\partial s^2} -(r-q)s\dfrac{\partial u}{\partial s}+ru=0, \quad (s,t)\in [0,\infty)\times[0,T].
\end{equation}
}
PDE \eqref{eq:ForwardBS_PDE} must be completed with initial and boundary conditions. The initial condition $u(s,0)$ depends on the payoff of the option.
%  and the boundary conditions should be carefully determined taking into account financial aspects as well as mathematical questions. Throughout the next subsections several types of options will be described, together with their corresponding initial and boundary conditions. 
% \subsection{Vanilla options}
% For example, a European call option is the right to buy a risky asset at a fixed strike price $K$ only at the future time $T$ (measured in years). The call option holder would exercise the option at expiry if the asset price is above the strike $K$ and not if it is below. Therefore, the payoff of a call option is $s_T-K$ if $s_T>K$ and $0$ otherwise. 
% Thus, the payoff of a European call option is $\max(s_T-K,0)$.
% Conversely, a put option gives the right to sell. At expiry the option is worth $\max(K-s_T,0)$. 

% The initial condition of \eqref{eq:ForwardBS_PDE} is $u(s,0) = \max(s-K,0)$ for call options and $u(s,0) = \max(K-s,0)$ for put options. 
% In order to solve numerically the Black-Scholes PDE we need to truncate the spatial domain. Therefore $u$ will be computed for $s\in(0,\bar{s})$, with $\bar{s}$ large enough. 
% % Besides, boundary conditions have to be imposed at the boundaries. 
% For call options the following Dirichlet boundary conditions can be used
% {
% \begin{align}
% &u(0,t)  = 0, \quad  &u(\bar{s},t)  &= \bar{s} e^{-q t} - Ke^{-rt},\, &\text{ for call options}
% \\
% &u(0,t) = Ke^{-rt} , \quad  &u(\bar{s},t)  &= 0,\,  &\text{ for put options}.
% \end{align}
% }
Standard European continuously monitored barrier options can be  priced in closed form. Their Greeks can be also computed analytically, see \cite{Merton73}. 
% See also \cite{Rubinstein-Reiner-91,Rich-94,Zhang98,Roman17}.
 Hereafter we are going to detail these formulas for down-and-in call options. Formulas for down-and-out call options can be inferred using that a portfolio consisting of an in option and its corresponding out option has the same price and Greeks of the corresponding vanilla option, i.e, $C(s,K,t)=C_{DO}(s,K,t) + C_{DI}(s,K,t)$. 
%  All these formulas are needed in order to measure the accuracy and the order of convergence of the proposed numerical schemes. Greek formulas are carefully detailed below since we were not able to find them in the literature.
Let $\bar{K}=\max(B,K)$ and  $\lambda = \frac{2}{\sigma^2}(r-q-\frac{\sigma^2}{2})$. The price of the down-and-in call can be computed in terms of the vanilla call/put prices:
{\small
\begin{eqnarray}
      C_{DI}(s,K,t) & =& \left(\dfrac{B}{s}\right)^{\lambda} \left[ C\left(\dfrac{B^2}{s},\bar{K}, t\right) + (\bar{K}-K) N\left(d_1\left( \dfrac{B^2}{s},\bar{K}\right) \right) \right]\nonumber \\
    & + &\left[ P(s,K,t)-P(s,B,t) + \frac{(B-K)\e{-rt}}{\sigma s \sqrt{t}}   N[-d_{1}(s,B)]\right] \mathds{1}_{B>K}.\label{eq:VDIC}
\end{eqnarray}
}
Where the analytical prices and Greeks of European vanilla call/put options are (see \cite{BS73,Merton73})
{\small
\begin{eqnarray}
C(s, K, t) &=& s e^{-qt}  N(d_1(s,K))-Ke^{-rt} N(d_2(s,K)), \label{BS-Callprice} \\
P(s, K, t) &=& K  e^{-rt} N(-d_2(s,K))-s e^{-qt} N(-d_1(s,K)),\label{BS-Putprice}
\end{eqnarray}
}
with $N$ the cumulative distribution of the standard normal, and $d_1$, $d_2$ given by:
{\small
$$
d_1(s,K)=\dfrac{1}{\sigma\sqrt{t}}\left[\ln\left(\dfrac{s}{K}\right)+\nu t \right],\quad  \nu =r-q+\dfrac{\sigma^2}{2}\label{d1},\quad d_2(s,K) =d_1(s,K)-\sigma\sqrt{t}.\label{d2} 
$$
}
% Hereafter we compute the delta and the gamma Greeks for the down-and-in call option. In the following expressions, for sake of brevity, in the formulas of the prices and deltas of vanilla call and put options, the time $t$ dependency is omitted. 
The Delta and the Gamma of the down-and-in call  can be computed by deriving \eqref{eq:VDIC} with respect to $s$, one or two times, respectively.
The Delta and the Gamma of the down-and-out call option can be obtained as $\Delta_{DO} = \Delta_C - \Delta_{DI}$ and $\Gamma_{DO} = \Gamma_C -\Gamma_{DI}.$
% \subsubsection{Asian options\label{sec:Asian}}
% Asian options are path dependent options whose payoff depends on the price $s_T$ of the risky asset and also on the arithmetic average price $a_T$ of the price $s_t$ defined by $a_t = \frac1t\int_0^t s_\tau d\tau$.
% Different types of Asian options are traded in financial markets. Floating strike call options have the payoff function $\max(s_T-a_T,0)$, while fixed strike call options consider the payoff $\max(a_T-K,0)$, $K$ being the strike price. American-style Asian options are also negotiated. 
% Let us denote by $u(s,a,t)$ the price of an Asian option. Under the standard Black-Scholes model for the risky asset, one can check that the price of an Asian option with payoff function $u_0(s,a)$ is the solution of the following forward in time two dimensional PDE (see \cite{Zvan-Forsyth-97})
% \begin{equation} \label{AsianPDE}
%     \dfrac{\partial u}{\partial t}  - \dfrac{1}{2}\sigma^2 s^2\dfrac{\partial^2 u}{\partial s^2} -rs\dfrac{\partial u}{\partial s}-\frac{1}{T-t}(s-a)\dfrac{\partial u}{\partial a}+ru=0, \quad u(s,a,0)=u_0(s,a).
% \end{equation}
% As an example, $u_0(s,a)=\max(a-K,0)$ is the initial condition for an European fixed strike call option.
A down-and-out call option can be priced solving PDE \eqref{eq:ForwardBS_PDE}  with initial condition
$$
u(s,0)=\begin{cases}
\max(s-K,0) & \text{for } s>B, \\
0 & \text{for } s\leq B,
\end{cases}
$$
in the localized domain $(s,t)\in[B,\bar{s}]\times(0,T]$ with the boundary conditions $u(B,t)=0$ and 
$u(\bar{s},t) = s \e{-q t} - K\e{-rt}$  for $t\in(0,T]$.

To solve this model numerically, first of all, the Black-Scholes PDE \eqref{eq:ForwardBS_PDE} is written in the conservative form  \eqref{ConsForm1D}, where the conservative functions are given by:
$$
f (u)=(\sigma^2-r+q)s u\, , \quad
g (u_s)=\dfrac{1}{2}\sigma^2s^2\dfrac{\partial u}{\partial s}\, , \quad
h(u)=(\sigma^2-2r+q)u.
$$
In this section a down-and-out call option with the market data $\sigma=0.2$, $r=0.05$, $q=0$, $T=1$, $K=70$ and the barrier at ${B=200}$ is priced. The computational domain is thus set to $[B,5B]$. 
In Figure \ref{fig:testBarrier} option prices, numerical errors, Deltas and Gammas are shown at $t=T$ considering a mesh with $N=800$. These plots show that the here proposed numerical methods are able to obtain good approximations without oscillations, even at difficult zones like close to the barrier. Table \ref{tb-Barrier} shows $L_1$ errors and $L_1$ order of convergence at $t=T$. Second order accuracy is achieved again. In this case, IMEX time step is between $200$ and $25606$ times larger than the explicit time step. Consequently, IMEX executes between $10$ and $12222$ times faster.

\vspace{-0.5cm}
\begin{figure}[!h]
    \centering \includegraphics[width=0.45\linewidth]{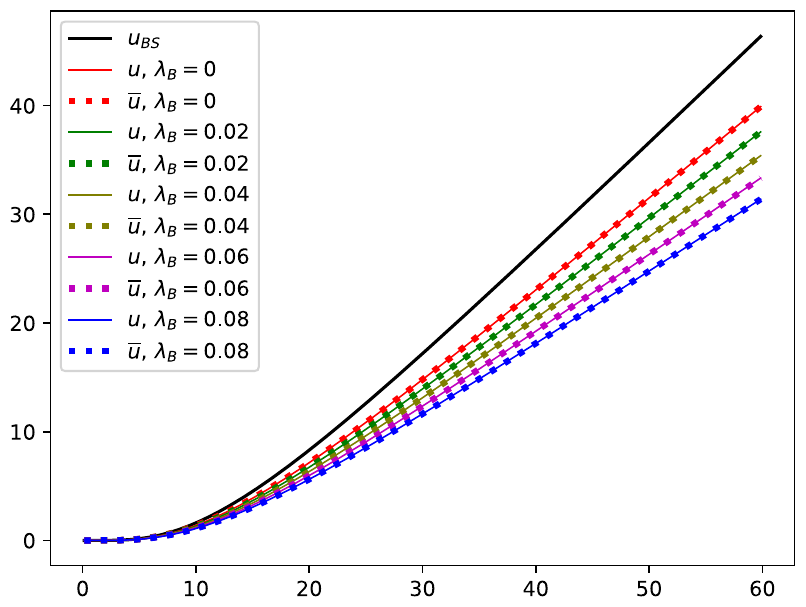}  
    \caption{Call option prices with valuation adjustments at $t=T$. Exact value in continous line, numerical value with dots, for different values of the default intensities $\lambda_B$}.
    \label{fig:testCVAprice}
    \end{figure}

\begin{figure}[!h]
    \vspace{-1cm}
\centering \includegraphics[width=0.45\linewidth]{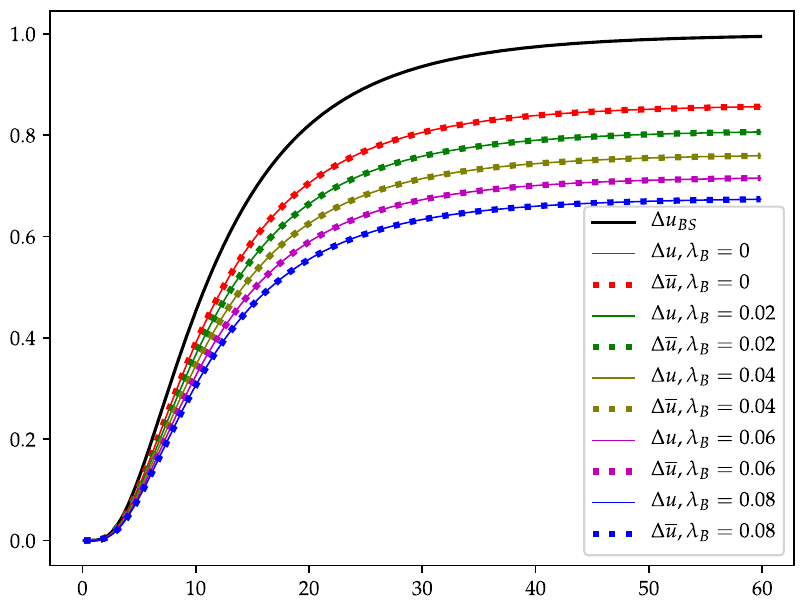} 
\centering \includegraphics[width=0.46\linewidth]{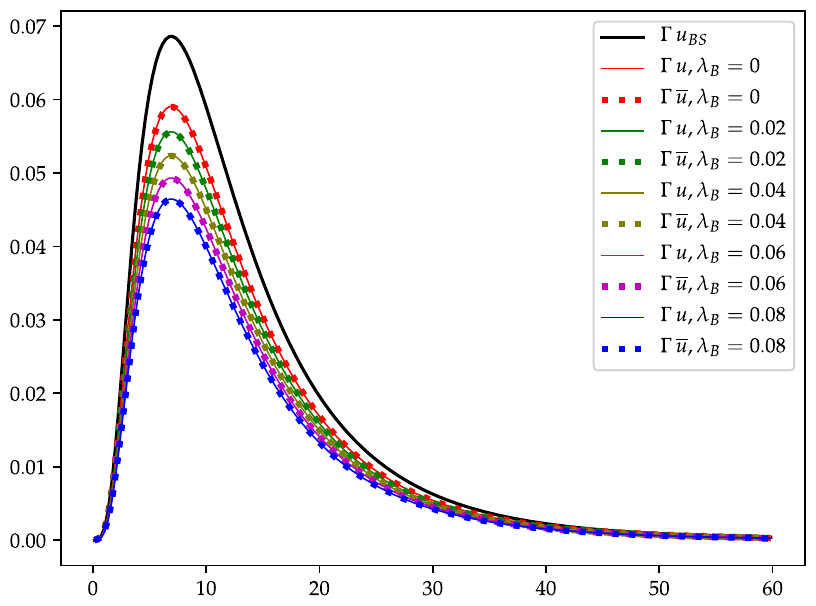}
\caption{$\Delta$, left, and $\Gamma$, right, with valuation adjustments at $t=T$. Exact value in continous line, numerical value with dots, for different values of the default intensities $\lambda_B$}
\label{fig:testCVAdelta}
\end{figure}

\vspace{-1cm}
\subsubsection{Nonlinear example. Counterparty credit risk}
\vspace{-0.05cm}
% So far we have been dealing with linear PDEs. 
In this section we present an important example of  a semi-linear PDE in finance, presented in \cite{BurgardKjaer}, that models the value of the derivative including counterparty risk:
\begin{align}
    &\dfrac{\partial u}{\partial t}  - \dfrac{1}{2}\sigma^2 s^2\dfrac{\partial^2 u}{\partial s^2} -(r-q)s\dfrac{\partial u}{\partial s}+ru=-(1-R_B)\lambda_B \min(u,0) \nonumber \\
    &\quad- ( (1-R_C) \lambda_C+ s_F) \max(u,0), \quad (s,t)\in[0,\infty)\times (0,T], \label{eq:ForwardBS_PDExva}
\end{align}
where $B$(buyer) and $C$(seller) are the two counterparties, $R_B$, $R_C \in [0,1]$ represent recovery rates on derivatives positions of parties $B$ and $C$, respectively, $\lambda_B$ and $\lambda_C$ are the default intensities and $s_F$ is the funding cost of the entity. Again analytic solution of the price and the Greeks can be computed and found in \cite{BurgardKjaer}.

Again Equation \eqref{eq:ForwardBS_PDExva} can be expressed in conservative form \eqref{ConsForm1D}, with the conservative functions given by:
\vspace{-0.25cm}
\begin{align*}
& f (u)=(\sigma^2-r+q)s u\, , \quad
g (u_s)=\dfrac{1}{2}\sigma^2s^2\dfrac{\partial u}{\partial s}\, , \quad \\
& h(u)=(\sigma^2-2r+q)u-(1-R_B)\lambda_B \min(u,0) - ( (1-R_C) \lambda_C+ s_F) \max(u,0).  
\end{align*}
The market data is taken as $T=5$, $K=15$, $r=0.02$, $\sigma=0.3$, $q=0$, $R_B=0.4$, $R_C=0.4$, $\lambda_B=\in\{0,0.02,0.04,0.06,0.08\}$, $\lambda_C=0.05$ and $s_F=(1-R_B)\lambda_B$. In Figures \ref{fig:testCVAprice}, \ref{fig:testCVAdelta} call option prices and Greeks are shown at $t=T$. Table \ref{tb-CVA} record $L_1$ errors and $L_1$ orders of convergence at $t=T$. The proposed numerical schemes continues to perform satisfactorily in this nonlinear setting: the second order is mantained and again a really accurate 
approximation of the price and the Greeks is obtained.

% \begin{figure}[!h]
% \centering \includegraphics[width=0.5\linewidth]{./images/counterparty/gamma.pdf}  
% \caption{$\Gamma$ with valuation adjustments at $t=T$. Exact value in continous line, numerical value with dots, for different values of the default intensities $\lambda_B$}
% \label{fig:testCVAgamma}
% \end{figure}

\vspace{-0.5cm}
\begin{table}[!h]
\begin{center}
{\scriptsize
\begin{tabular}{|c||c|c|c|c|}
\hline
& \multicolumn{4}{c|}{IMEX}  \\
\hline
$N$  & $L_1$ \text{ error} &  \text{ Order}& $\Delta t$ & Time (s) \\
\hline
$50$ &  $1.4323 \times 10^{-1}$ & $--$& $1.44\times 10^{-1}$ & $  7.2 \times 10^{-3}$ \\ 
$100$ & $3.6714 \times 10^{-2}$ & $1.96$& $7.17\times 10^{-2}$ & $ 1.8\times 10^{-3}$ \\
$200$ & $9.2457 \times 10^{-3}$ &  $1.99$ & $3.58\times 10^{-2}$ &  $  4.5\times 10^{-3}$\\
$400$ & $2.3140 \times 10^{-3}$ & $2.00$ &  $1.78\times 10^{-2}$& $ 1.1 \times 10^{-2}$ \\
$800$ &  $5.7768 \times 10^{-4}$ & $2.00$ & $8.93\times 10^{-3}$ & $  2.4\times 10^{-2}$ \\
$1600$ &  $1.4413 \times 10^{-4}$ & $2.00$ & $4.46 \times 10^{-3}$ & $ 7.8\times 10^{-2}$ \\
$3200$ &  $3.5943 \times 10^{-5}$ & $2.00$ & $2.23 \times 10^{-3}$ &  $ 2.7\times 10^{-1}$ \\
$6400$ &  $8.9052 \times 10^{-6}$ & $2.01$ & $1.11 \times 10^{-3}$ & $ 1.0\times 10^{0}$ \\
\hline
\end{tabular} 
\begin{tabular}{|c||c|c|c|c|}
    \hline
& \multicolumn{4}{c|}{Explicit}  \\
\hline
$N$  & $L_1$ \text{ error} &  \text{ Order}& $\Delta t$ & Time (s) \\
\hline
$50$ & $1.4255 \times 10^{-1}$ & $--$ & $2.23 \times 10^{-3}$ & $ 1.2\times 10^{-2}$ \\
$100$ & $3.5607 \times 10^{-2}$ & $2.00$ & $5.56 \times 10^{-4}$ & $ 1.7\times 10^{-2}$\\
$200$ &  $8.8734 \times 10^{-3}$ & $2.00$ & $1.39 \times 10^{-4}$ & $ 8.3 \times 10^{-1}$ \\
$400$ &  $2.2145 \times 10^{-3}$ & $2.00$ & $3.47 \times 10^{-5}$ & $ 3.8\times 10^{0}$\\
$800$ & $5.5312 \times 10^{-4}$ &$2.00$ & $8.68 \times 10^{-6}$ & $ 2.2\times 10^{1}$\\
$1600$ & $1.3823 \times 10^{-4}$ &$2.00$ & $2.17 \times 10^{-6}$ & $ 8.7\times 10^{1}$\\
$3200$ & $3.4467 \times 10^{-5}$ & $2.00$ & $5.42 \times 10^{-7}$  & $ 7.3\times 10^{2}$ \\
$6400$ & $8.6169 \times 10^{-6}$ & $2.00$ & $1.36 \times 10^{-7}$ & $ 3.5\times 10^{3}$ \\
\hline
\end{tabular} 
\caption{$L_1$ errors and $L_1$ orders of convergence of the IMEX and explicit finite volume methods for the call with valuation adjustments.}
\label{tb-CVA} 
}
\end{center}
\end{table}

\vspace{-1.5cm}
\section{Conclusions}
\vspace{-0.25cm}
In this article we have shown that finite volume IMEX Runge-Kutta numerical sche\-mes are very suitable for solving PDE option pricing  problems. On the one hand, the IMEX time discretization is remarkably efficient. Indeed, large time steps can be used, avoiding the need to use the smaller, and possibly extremely small, time steps enforced by the diffusion stability condition, which has to be satisfied in explicit schemes. Numerical results show that IMEX outperforms the explicit method. In fact, IMEX is the only way to solve problems in highly refined meshes is space. 
% Besides, even in its worst scenarios, IMEX performs at least as well as the explicit method.
On the other hand, finite volume space discretization allows to obtain second order convergence in convection dominated problems and/or problems with non smooth initial and/or boundary conditions, which is the usual situation in finance.  Thus, no special regularization techniques of the non smooth data need to be taken into account. The accuracy of the numerical scheme turns to be of key importance for the accurate and non oscillatory computation of the Greeks. 
% Finally, in this paper we provide a set of benchmark problems, together with their analytical solutions. These benchmarks can also be valuable for mathematical researchers working in the development of high order numerical schemes for advection-diffusion problems.

\section*{Acknowledgements}

 % been funded by FEDER and the Spanish Government  through the coordinated Research project RTI2018-096064-B-C1; by the Junta de Andaluc\'ia research projects P18-RT-3163 and the Junta de Andalucia-FEDER-University of M\'alaga research project UMA18-FEDERJA-16; and the University of M\'alaga. 

M. Castro has has been partially supported by the grant PDC2022-133663-C21 funded by MCIN/AEI/10.13039/501100011033 and “European Union NextGenerationEU/PRTR"  and  the grant \mbox{PID2022-137637NB-C21} funded by MCIN/AEI/\-10.13039/50110001103 and “ERDF A way of making Europe”.
The other authors' research has been funded by the Spanish MINECO under research project number PDI2019-108584RB-I00 and by the grant ED431G 2019/01 of CITIC, funded by Conseller\'ia de Educaci\'on, Universidade e Formaci\'on Profesional of Xunta de Galicia and FEDER.

\end{document}